\documentclass[number,citesort,seceqn,dvips]{arxbj}

\aid{0}
\volume{18}
\issue{1}
\pubyear{2012}
\firstpage{24}
\lastpage{45}
\doi{10.3150/10-BEJ327}

\makeatletter
\newproclaim{definition}{Definition}[section]
\newtheorem{theorem}{Theorem}[section]
\newremark{examples}{Examples}
\newtheorem{proposition}[theorem]{Proposition}
\newtheorem{lemma}[theorem]{Lemma}

\newcommand{\1}{{\mathbf1}}
\newcommand{\la}{\lambda}

\newcommand{\cc}{\mathcal C}
\newcommand{\cf}{\mathcal F}
\newcommand{\cl}{\mathcal L}

\newcommand{\al}{\alpha}
\newcommand{\si}{\sigma}

\newcommand{\N}{\mathbb N}
\newcommand{\R}{\mathbb R}
\newcommand{\G}{\mathbb P}

\makeatother

\begin{document}
\begin{frontmatter}

\title{Stochastic delay equations with non-negativity constraints
driven by fractional Brownian motion}

\runtitle{Stochastic delay equations with constraints}

\begin{aug}
\author{\fnms{Mireia} \snm{Besal\'{u}}\thanksref{e1}\ead[label=e1,mark]{mbesalu@ub.edu}} \and
\author{\fnms{Carles} \snm{Rovira}\thanksref{e2}\corref{}\ead[label=e2,mark]{Carles.Rovira@ub.edu}}
\runauthor{M. Besal\'{u} and C. Rovira}
\address{Facultat de Matem\`{a}tiques, Universitat de Barcelona, Gran
Via 585, 08007 Barcelona, Spain.\break\printead{e1,e2}}
\end{aug}

\received{\smonth{3} \syear{2010}}
\revised{\smonth{8} \syear{2010}}

%
\begin{abstract}
In this note we prove an existence and uniqueness result for the
solution of multidimensional stochastic delay differential equations
with normal reflection. The equations are driven by a~fractional
Brownian motion with Hurst parameter $H > 1/2$. The stochastic integral
with respect to the fractional Brownian motion is a pathwise
Riemann--Stieltjes integral.
\end{abstract}

%

\end{frontmatter}
%

\section{Introduction}

Consider a stochastic delay differential equation with positivity
constraints. More precisely, we deal with a stochastic delay
differential equation with normal reflection on $\R^d$ of the form
\begin{eqnarray}\label{equation}
X(t)&=&\eta(0)+\int_0^t b(s,X)\,\mathrm{d}s+\int_0^t \si
\bigl(s,X(s-r)\bigr)\,\mathrm{d}W_s^H+Y(t),\qquad t\in(0,T],\nonumber\\ [-8pt]\\ [-8pt]
X(t)&=& \eta(t),\qquad t\in[-r,0].\nonumber
\end{eqnarray}
Here, $r$ denotes a strictly positive time delay,
$W^H=\{W^{H,j},j=1,\ldots,m\}$ are independent fractional
Brownian motions with Hurst parameter $H>\frac{1}{2}$ defined in a
complete probability space $(\Omega,\cf,\G)$, the hereditary term
$b(s,X)$ depends on the path $\{X(u),-r\leq u\leq s\}$,
while $\eta\dvtx[-r,0]\rightarrow\R^d_+$ is a non-negative smooth function
with $\R^d_+=\{u\in\R^d;u_i\geq
0 \mbox{ for } i=1,\ldots,d\}$ and $Y$ is a vector-valued
non-decreasing process which ensures that the non-negativity
constraints on X are enforced.

Set
\begin{equation}\label{num1}
Z(t)= \eta(0)+\int_0^t b(s,X)\,\mathrm{d}s+\int_0^t \si\bigl(s,X(s-r)\bigr)\,\mathrm{d}W_s^H,\qquad
t\in[0,T].
\end{equation}
It is known that we have an explicit formula for the regulator term $Y$
in terms of $Z$: for each $i=1,\ldots,d,$
\[
Y^i(t)=\max_{s\in[0,t]} (Z^i(s))^-,\qquad t\in[0,T].
\]
The solution of (\ref{equation}) then satisfies
\[
X(t)=
\cases{
Z(t)+Y(t),&\quad $t\in[0,T]$,\cr
\eta(t),&\quad $t\in[-r,0]$.
}
\]

We call (\ref{equation}) a delay differential equation with reflection
with hereditary drift driven by a~fractional Brownian motion. To the
best of our knowledge, this problem has not been considered before in
the extensive literature on stochastic differential
equations.

There are many references to stochastic systems with delay (see, as a
basic reference,~\cite{M}), but the literature concerning stochastic
differential equations with delay driven by a fractional Brownian
motion is scarce. The existence and uniqueness of solutions \cite{FR,LT,NNT,TT}, existence and regularity of the
density \cite{LT} and convergence when the delay goes to zero
\cite{FerRo} have all been studied.

On the other hand, there has been little work on stochastic
differential equations with delay and non-negativity constraints. We
can only refer the reader to the book by Kushner \cite{K} dedicated to
the study of numerical methods for this class of equations, and the
paper of Kinnally and Williams \cite{KW}, where the authors obtain
sufficient conditions for existence and uniqueness of stationary
solutions for stochastic differential equations with delay and
non-negativity constraints driven by standard Brownian motion.

As described in the paper of Kinnally and Williams \cite{KW}, there
are some models affected by some types of noise where the dynamics are
related to propagation delay and some of them are naturally
non-negative quantities, for instance, applications such as rates and
prices in Internet models, and concentrations of ions or proportions of
a population that are infected (see the references in \cite{KW}).
Therefore, it is natural to consider stochastic differential equations
with delay and non-negativity constraints. In this paper, we initiate
the study when the noise is not a standard Brownian motion but a
fractional Brownian motion.

The main novelty of this paper is the use of non-negative constraints
dealing with fractional Brownian motion. We have used
Skorokhod's mapping.

 Set
\[
\cc_+(\R_+,\R^d):=\{x\in\cc(\R_+,\R^d)\dvtx
x(0)\in\R^d_+\}.
\]
We now recall the
Skorokhod problem.

\begin{definition}
Given a path $z\in\cc_+(\R_+,\R^d)$, we say that a pair $(x,y)$ of
functions in $\cc_+(\R_+,\R^d)$ solves the Skorokhod problem for $z$
with reflection if:
\begin{enumerate}
\item $x(t)=z(t)+y(t)$ for all $t\geq0$ and $x(t)\in\R_+^d$ for each
$t\geq0$;
\item for each $i=1,\ldots, d$, $y^i(0)=0$ and $y^i$ is
non-decreasing;
\item for each $i=1,\ldots, d$, ${ \int_0^t
x^i(s)\,\mathrm{d}y^i(s)=0}$ for all $t\geq0$, so $y^i$ can increase only when~$x^i$ is at zero.
\end{enumerate}
\end{definition}

 It is known that we have an explicit formula for $y$ in terms
of $z$: for each $i=1,\ldots,d,$
\[
y^i(t)=\max_{s\in[0,t]} (z^i(s))^-.
\]

The path $z$ is called the \textit{reflector} of $x$ and the path $y$ is
called the \textit{regulator} of $x$. We use the Skorokhod mapping for
constraining a continuous real-valued function to be non-negative by
means of reflection at the origin. We will apply it to each path of~$z$
defined by (\ref{num1}). Note that because we are dealing with a
multidimensional case, the mapping will be applied to each component.

We must also explain how to understand the stochastic integral
appearing in (\ref{equation}). Since $H > \frac12$, the stochastic
integral in (\ref{equation}) is defined using a pathwise approach.
Indeed, if we have a~stochastic process $\{u(t), t \ge0\}$ whose
trajectories are $\lambda$-H\"{o}lder continuous with $\lambda
> 1-H$, then the Riemann--Stieltjes integral $\int_0^T u(s)\,\mathrm{d}W_s^H$
exists for each trajectory (see Young \cite{Y}). Using the techniques
introduced by Young \cite{Y} and the $p$-variation norm, Lyons \cite{L}
began the study of integral equations driven by functions with bounded
$p$-variation, where $p\in[1,2)$. Z\"{a}hle \cite{Z} introduced a
generalized Stieltjes integral using the techniques of fractional
calculus. The integral is expressed in terms of fractional derivative
operators and coincides with the Riemann--Stieltjes integral $\int_0^T
f\,\mathrm{d}g$ when the functions $f$ and $g$ are H\"{o}lder continuous of
orders $\lambda$ and $\beta$,\vspace*{1pt} respectively, with $\lambda+ \beta> 1$.
Using this Riemann--Stieltjes integral, Nualart and Rascanu~\cite{NR}
obtained the existence and uniqueness of a solution for a class of
multidimensional integral equations.

In this paper, also using the Riemann--Stieltjes integral, we will
prove the existence and uniqueness of a solution to
equation (\ref{equation}). Our results are inspired by those in
Nualart and Rascanu \cite{NR} and Ferrante and Rovira \cite{FerRo}.
Using some estimates presented in those papers, we will first prove our
results for deterministic equations and will then easily apply them
pathwise to fractional Brownian motion.

Since our definition of the stochastic integral holds for
$H>\frac{1}{2}$, we cannot extend our approach to the case
$H\in(\frac{1}{3},\frac{1}{2})$. However, in a forthcoming
paper, we will use the method used by Hu and Nualart in \cite{HN} to
consider the case $H<\frac{1}{2}$.

The structure of the paper is as follows. In the next section, we give
our hypothesis and we state the main results of our paper. In Section
\ref{sec3}
we give some useful estimates for Lebesgue
and Riemann--Stieltjes integrals inspired by the results in \cite{NR}
and~\cite{FerRo}.
Section~\ref{sec4} is devoted to proving our main result: the existence,
uniqueness and boundedness of solutions to deterministic equations. In
Section \ref{sec5} we recall how to apply the deterministic results to the
stochastic case, while the \hyperref[appendix]{Appendix} is devoted to giving some technical
results such as a fixed point theorem and some properties related to
the Skorokhod problem.\looseness=-1

\section{Main results}\label{sec2}

Let $\alpha\in(0,\frac{1}{2})$ and $r>0$. Let
$(s,t)\subseteq[-r,T]$ and denote by $W_0^{\alpha,\infty}(s,t;\R^d)$
the space of measurable functions $f\dvtx[s,t]\rightarrow\R^d$ such that
\[
\|f\|_{\alpha,\infty(s,t)}:=\sup_{u\in[s,t]}
\biggl(|f(u)|+\int_{s}^u\frac{|f(u)-f(v)
|}{(u-v)^{\alpha+1}}\,\mathrm{d}v\biggr)<\infty.
\]
For any $0<\lambda\leq1$, denote by $C^{\lambda}(s,t;\R^d)$ the space
of $\lambda$-H\"{o}lder continuous functions~$f\dvtx\allowbreak[s,t]\rightarrow\R^d$
such that
\[
\|f\|_{\lambda(s,t)}:=\|f\|_{\infty
(s,t)}+\sup_{s\leq u<v\leq t}\frac{|f(v)-f(u)
|}{(v-u)^{\lambda}}<\infty,
\]
where
\[
\|f\|_{\infty(s,t)}:=\sup_{u\in[s,t]}|f(u)|.
\]
Since we will use the spaces $W_0^{\alpha,\infty}(-r,T;\R^d)$ and
$C^\lambda(-r,T;\R^d)$ extensively, we will use the notation
$\|f\|_{\alpha,\infty(r)}:=\|f\|_{\alpha
,\infty(-r,T)},
\|f\|_{\alpha,\lambda(r)}:=\|f\|_{\alpha
,\lambda(-r,T)},
\|f\|_{\lambda(r)}:= \|f\|_{\lambda(-r,T)}$ and
$\|f\|_{\infty(r)}:=\|f\|_{\infty(-r,T)}.$

Consider the following hypothesis.
\begin{itemize}
\item[(H1)] $\sigma\dvtx[0,T]\times\R^d\rightarrow
\R^d\times\R^m$ is a measurable function such that there exist some
constants $\beta>0$ and $M_0>0$ such that the following properties
hold:
\begin{enumerate}
\item $|\sigma(t,x)-\sigma(t,y)|\leq
M_0|x-y|\ \forall x,y\in\R^d, \forall t\in[0,T]$;
\item
$|\sigma(t,x)-\sigma(s,x)|\leq
M_0|t-s|^\beta\ \forall x\in\R^d, \forall t,s\in[0,T]$.
\end{enumerate}
\item[(H2)] $b\dvtx[0,T]\times C(-r,T;\R^d)\rightarrow\R^d$
is a
measurable function such that for every $t>0$ and $f\in C(-r,T;\R^d),
b(t,f)$ depends only on $\{f(s);-r\leq s\leq t\}$.
Moreover, there exists some $b_0\in L^\rho(0,t;\R^d)$ with $\rho\geq2$
and $\forall N\geq0$ there exists some $L_N>0$ such that:
\begin{enumerate}
\item$|b(t,x)-b(t,y)|\leq L_N\sup_{-r\leq s\leq
t}|x(s)-y(s)|\ \forall x,y \mbox{ such that }
\|x\|_{\infty(r)}\leq N,\allowbreak\|y\|_{\infty
(r)}\leq
N\  \forall t\in[0,T]$;
\item$|b(t,x)|\leq L_0\sup_{-r\leq
s\leq t}|x(s)|+b_0(t)\ \forall t\in[0,T]$.
\end{enumerate}
\item[(H3)] There exist some $\gamma\in[0,1]$ and $K_0>0$
such that
\[
|\si(t,x)|\leq K_0 (1+|x|^\gamma)\qquad\forall x\in\R^d,\ \forall
t\in[0,T].
\]
\end{itemize}
Under these assumptions, we are able to prove that our problem admits a
unique solution. Our main result reads as follows.
\begin{theorem}
Assume that $\eta\in W_0^{\alpha,\infty}(-r,0;\R^d_+),$ and that
$b$ and $\sigma$ satisfy hypotheses~\textup{(H1)} and~\textup{(H2)}, respectively, with $\beta>1-H$. Set $\alpha
_0:=\min\{\frac{1}{2}, \beta\}$.
If $\alpha\in(1-H,\alpha_0)$ and $\rho\leq\frac{1}{\alpha}$,
then the
equation (\ref{equation}) has a unique solution,
\[
X\in L^0(\Omega,\cf,\G;W_0^{\alpha,\infty}(-r,T;\R^d)),
\]
and for $P$-almost all $\omega\in\Omega,$ $X(\omega,\cdot) \in
C^{1-\al}(0,T;\R^d)$.

Moreover, if $\alpha\in(1-H, \al_0 \vee(2-\gamma)/4)$ and \textup{(H3)}
holds, then $E ( \Vert X \Vert^p_{\al,\infty(r)} ) < \infty\
\forall p \ge1.$\
\end{theorem}

\begin{examples*} Note that the following equations
satisfy our hypothesis:
\begin{longlist}
\item[(a)] (linear example) for any $a, b \in\R,$
\begin{eqnarray*}
X(t)&=&r +\int_0^t X(s-r)\,\mathrm{d}s+\int_0^t \bigl(a X(s-r) + b\bigr)\,\mathrm{d}W_s^H+Y(t),\qquad t\in(0,T],\\
X(t)&=& t+r,\qquad  t\in[-r,0];
\end{eqnarray*}
\item[(b)] (non-linear example)
\begin{eqnarray*}
X(t)&=&\int_0^t \cos(X(s))\,\mathrm{d}s+\int_0^t \sin\bigl(s+X(s-r)\bigr)\,\mathrm{d}W_s^H+Y(t),\qquad t\in(0,T],\\
X(t)&=& t^2,\qquad  t\in[-r,0].
\end{eqnarray*}
\end{longlist}
\end{examples*}

\section{Preliminaries}\label{sec3}

In this section, we give some useful estimates for Lebesgue and
Riemann--Stieltjes integrals. These types of estimates were presented
in the work of Nualart and Rascanu \cite{NR} and adapted to the delay
case by Ferrante and Rovira \cite{FerRo}. Since our results are directly
inspired by these works, we do not give the proofs, but instead direct
the reader to these references.

We will need to
introduce a new norm in the space $W_0^{\alpha,\infty}(s,t;\R^d)$: for
any $\lambda\geq1$,
\[
\|f\|_{\alpha,\lambda(s,t)}:=\sup_{u\in[s,t]}\exp
(-\lambda u)\biggl(|f(u)|+\int_{s}^u \frac{
|f(u)-f(v)|}{(u-v)^{\alpha+1}}\,\mathrm{d}v\biggr).
\]
It is easy to check that for any $\lambda\geq1$, this norm is
equivalent to $\|f\|_{\alpha,\infty(s,t)}$.

\subsection{Lebesgue integral}

We first consider the ordinary Lebesgue integral. Given a measurable
function $f\dvtx[-r,T]\rightarrow\R^d$, we define
\[
F^{(b)}(f)(t)=\int_0^t b(u,f)\,\mathrm{d}u.
\]
We first recall some estimates that constitute an obvious adaptation
of \cite{FerRo}, Proposition~2.2.
\begin{proposition} \label{Prop2.2}
Assume that b satisfies \textup{(H2)} with
$\rho=\frac{1}{\alpha}$ and $[s,t]\subseteq[0,T]$. If $f\in
W_0^{\alpha,\infty}(-r, t;\R^d),$ then $F^{(b)}(f)(\cdot)=\int_0^.
b(u,f)\,\mathrm{d}u\in C^{1-\alpha}(s,t;\R^d)$ and\\
\hspace*{11pt} 1. ${ \bigl\|F^{(b)}(f)\bigr\|_{1-\alpha(s,t)}\leq
d^{(1)}\bigl(1+\|f\|_{\infty(-r,t)}\bigr),}$\\
\hspace*{11pt} 2. $\displaystyle \bigl\|F^{(b)}(f)\bigr\|_{\alpha,\lambda(s,t)}\leq
d^{(2)}\biggl(\frac{1}{\lambda^{1-2\alpha}}+\frac{\|f\|
_{\alpha,\lambda(-r,t)}}{\lambda^{1-\alpha}}\biggr)$\\
for all $\lambda\geq1$, where $d^{(i)}$, $i\in\{1,2\}$, are
positive constants depending only on $\alpha,t,L_0$ and
$B_{0,\alpha}=\|b_0\|_{L^{1/\alpha}}$.
\end{proposition}
%

\subsection{Riemann--Stieltjes integral}\vspace*{3pt}

Let us now consider the Riemann--Stieltjes integral introduced by
Z\"{a}hle, which is based on fractional integrals and derivatives. We
refer the reader to the paper of Z\"{a}hle \cite{Z} and the references
therein for a detailed presentation of this generalized Stieltjes
integral and the associated fractional calculus. Here, we will just
recall some basic facts.

Fix a parameter
$0<\al<\frac{1}{2}$. Denote by $W_T^{1-\al,\infty}(0,T;\R)$ the space
of measurable functions $g\dvtx[0,T]\rightarrow\R$ such that
\[
\|g\|_{1-\al,\infty,T}:=\sup_{0<s<t<T}\biggl(\frac
{|g(t)-g(s)|}{(t-s)^{1-\al}}+\int_s^t\frac{
|g(y)-g(s)|}{(y-s)^{2-\al}}\,\mathrm{d}y\biggr)<\infty.
\]
Moreover, if $g$ belongs to $W_T^{1-\al,\infty}(0,T;\R)$, then we
define
\begin{eqnarray*}
\Lambda_\al(g)&:=&\frac{1}{\Gamma(1-\al)}\sup_{0<s<t<T}
|(D^{1-\al}_{t^-}
g_{t^-})(s)|\\
&\leq&\frac{1}{\Gamma(1-\al)\Gamma(\al)}\|g\|_{1-\al,\infty,T}<\infty,
\end{eqnarray*}
where $D^{1-\alpha}_{t^-} g_{t^-}$ denotes a fractional Weyl
derivative,
\[
g_{t-}(s)=\bigl(g(s)-g(t-)\bigr)\1_{(0,t)}(s),
\]
where
$g(t-)=\lim_{\varepsilon\searrow0} g(t-\varepsilon)$ and $\Gamma$ is
the Euler function. We also denote by~$W_0^{\al,1}(0,\allowbreak T;\R)$ the space
of measurable functions $f$ on $[0,T]$ such that
\[
\|f\|_{\al,1}:=\int_0^T\frac{|f(s)|}{s^\al
}\,\mathrm{d}s+\int_0^T\int_0^s \frac{|f(s)-f(y)|}{(s-y)^{\al
+1}}\,\mathrm{d}u\,\mathrm{d}s<\infty.
\]
Note that if $f$ is a function in the space $W_0^{\al,1}(0,T;\R)$ and
$g\in W_T^{1-\al,\infty}(0,T;\R),$ then the integral $\int_0^t f\,\mathrm{d}g$
exists for all $t\in[0,T]$ and we can define
\[
G(f)(t):=\int_0^t f(s)\,\mathrm{d}g_s=\int_0^T f(s) \1_{(0,t)}(s)\,\mathrm{d}g_s.
\]
Furthermore, the following estimate holds:
\begin{equation}\label{cotasi1}
\biggl|\int_0^t f\,\mathrm{d}g\biggr|\leq\Lambda_\al(g)\|f\|
_{\al,1}.
\end{equation}
Moreover, if $f \in W_0^{\alpha,\infty} (0,T)$, it is proved in
\cite{NR}, Proposition 4.1, that for each $s<t$,
\begin{equation}\label{cotasi2}
\biggl|\int_s^t f\,\mathrm{d}g\biggr|\leq\Lambda_\al(g)C_{\al,T}
(t-s)^{1-\al}\|f\|_{\al,\infty}.
\end{equation}
Let us consider the term
\[
G_r^{(\sigma)}(f)(t)=\int_0^t \sigma\bigl(s,f(s-r)\bigr)\,\mathrm{d}g_s.
\]
For the Riemann--Stieltjes integral, we will also give a version of
\cite{FerRo}, Proposition 2.4.
\begin{proposition} \label{Prop2.4}
Let $g\in W_T^{1-\al,\infty}(0,T)$. Assume that $\sigma$ satisfies
\textup{(H1)} and $[s,t]\subseteq[0,T]$. If $f\in
W_0^{\alpha,\infty}(-r,T;\R^d),$ then
\[
G_r^{(\sigma)}(f)\in C^{1-\alpha}(s,t;\R^d)\subset W_0^{\alpha
,\infty}(s,t;\R^d)
\]
and\\
\hspace*{11pt}1. $ {\bigl\|G_r^{(\sigma)}(f)\bigr\|_{1-\alpha
(s,t)}\leq
\Lambda_\alpha(g)d^{(3)}\bigl(1+\|f\|_{\alpha,\infty(-r,t-r)}\bigr)}$,\\
\hspace*{11pt}2.$ {\bigl\|G_r^{(\sigma)}(f)\bigr\|_{\alpha,\lambda
(s,t)}\leq
\frac{\Lambda_\alpha(g)d^{(4)}}{\lambda^{1-2\alpha}}\bigl(1+\|
f\|_{\alpha,\lambda(-r,t-r)}\bigr)}$\\
for all $\lambda\geq1$, where $d^{(i)}$, $i\in\{3,4\}$, are
positive constants independent of $\lambda$, f and g.
\end{proposition}

Finally, we recall \cite{FerRo}, Proposition 2.6.
Consider $\varphi(\gamma,\al)$ defined such that
$\varphi(\gamma,\al)=2\al$ if $\gamma=1$,
$\varphi(\gamma,\al)>1+\frac{2\al-1}{\gamma}$ if
$\frac{1-2\al}{1-\al}\leq\gamma<1$ and $\varphi(\gamma,\al)=\al
$ if
$0\leq\gamma<\frac{1-2\al}{1-\al}$. Note that
$\varphi(\gamma,\al)\in[\al,2\al]$.

\begin{proposition} \label{Prop2.6}
Let $g\in W_T^{1-\al,\infty}(0,T)$. Assume that $\sigma$ satisfies
\textup{(H1)} and \textup{(H3)}. If $f\in
W_0^{\alpha,\infty}(-r,T;\R^d),$ then
\[
\bigl\|G_r^{(\si)}(f)\bigr\|_{\al,\la}\leq\Lambda_\al
(g)d^{(5)}\biggl(1+\frac{\|f\|_{\al,\la(r)}}{\la
^{1-\varphi(\gamma,\al)}}\biggr)
\]
for all $\la\geq1$, where $d^{(5)}$ is a positive constant depending
only on $\al,\beta,T,d,m$ and
$B_{0,\al}=\|b_0\|_{L^{1/\al}}$.
\end{proposition}

\vspace*{-3pt}\section{Deterministic integral equations}\vspace*{-3pt}\label{sec4}

In this section, we give all the deterministic results.

For simplicity, let us assume that $T=Mr$. Set
$0<\alpha<\frac{1}{2},g\in W_T^{1-\alpha,\infty}(0,T;\R^d)$ and
$\eta\in W_0^{\alpha,\infty}(-r,0;\R^d_+)$. Consider the deterministic
stochastic differential equation on $\R^d$
\begin{eqnarray}\label{det}
x(t)&=&\eta(0)+\int_0^t b(s,x)\,\mathrm{d}s+\int_0^t\sigma\bigl(s,x(s-r)\bigr)\,\mathrm{d}g_s+
y(t),\qquad t\in(0,T],\nonumber\\ [-6pt]\\ [-6pt]
x(t) &=& \eta(t),\qquad t\in[-r,0],\nonumber
\end{eqnarray}
where, for each $i=1,\ldots,d,$
\[
{\mathit y^i(t)=\max_{s\in[0,t]} (z^i(s))^-,\qquad
t\in[0,T],}
\]
and
\[
{\mathit z(t)= \eta(0)+\int_0^t b(s,x)\,\mathrm{d}s+\int_0^t
\si\bigl(s,x(s-r)\bigr)\,\mathrm{d}g_s,\qquad t\in[0,T].}
\]
The existence and uniqueness result reads as follows.
\begin{theorem}\label{teoexi}
Assume that $b$ and $\sigma$ satisfy hypotheses
\textup{(H1)} and \textup{(H2)}, respectively, with $\rho=1/\alpha$
and $0 < \alpha<\min\{\frac12, \beta\}$. The equation (\ref{det})
then has a unique solution $x\in W_0^{\alpha,\infty}(-r,T;\R^d_+)$.
\end{theorem}
\begin{pf}
To prove that equation (\ref{det}) admits a unique solution on
$[-r,T]$, we shall use an induction argument. We will prove that if
equation (\ref{det}) admits a unique solution on $[-r,nr]$, then we
can further prove that there is a unique solution on the interval
$[-r,(n+1)r]$.

Our induction hypothesis, for $k\leq M,$ is the following:\vspace*{11.4pt}

($\mathrm{H_k}$) \textit{The equation
\begin{eqnarray*}
 x^k(t)&=&\eta(0)+\int_0^t b(s,x^k)\,\mathrm{d}s+\int_0^t\sigma
\bigl(s,x^k(s-r)\bigr)\,\mathrm{d}g_s+y_k(t),\qquad
 t\in[0,kr],\\
 x^k(t) &=& \eta(t),\qquad t\in[-r,0],
\end{eqnarray*}
where, for each $i=1,\ldots,d,$
\[
y^i_k(t)=\max_{s\in[0,t]} (z^i_k(s))^-,\qquad
t\in[0,kr],
\]
with
\[
 z_k(t)= \eta(0)+\int_0^t b(s,x^{k})\,\mathrm{d}s+\int_0^t
\si\bigl(s,x^k(s-r)\bigr)\,\mathrm{d}g_s,\qquad t\in[0,kr],
\]
has a unique solution $x^k\in
W_0^{\alpha,\infty}(-r,kr;\R^d_+)$.}\vspace*{11.4pt}

The initial
case can be easily checked. Assume now that ($\mathrm{H_i}$) is true for
all $i\leq n$, where $n<M$. We wish to check ($\mathrm{H_{n+1}}$).

Clearly, for $t \in[-r,nr]$, $x^{n+1}(t)$ will coincide
with $x^n(t)$, the solution of the equation of ($\mathrm{H_n}$). Moreover,
for $t \in[-r,nr]$, $y_{n+1}(t)$ will coincide with $y_n(t)$.
We can therefore write the equation of ($\mathrm{H_{n+1}}$) as\vspace*{-1pt}
\begin{eqnarray}\label{equationn}
x^{n+1}(t)&=&\eta(0)+\int_0^t b(s,x^{n+1})\,\mathrm{d}s+\int_0^t\sigma
\bigl(s,x^n(s-r)\bigr)\,\mathrm{d}g_s\nonumber\\
&&{}+y_{n+1}(t),\qquad t\in[0,(n+1)r],\\
x^{n+1}(t) &=& \eta(t),\qquad t\in[-r,0].\nonumber
\end{eqnarray}
Moreover, using the notation introduced in the previous
section, we have
\begin{eqnarray*}
x^{n+1}(t)&=&\eta(0)+F^{(b)}(x^{n+1})+G^{(\si)}(x^n)+y_{n+1}(t),\qquad
t\in[0,(n+1)r],\\
x^{n+1}(t) &=& \eta(t),\qquad t\in[-r,0].
\end{eqnarray*}
%

 The proof will be divided into three steps:
\begin{enumerate}
\item if $x^{n+1}$ is a solution of $\mathrm{(H_{n+1})}$ in the space
$\cc(-r,(n+1)r;\R^d_+),$ then $x^{n+1}\in
W_0^{\al,\infty}(-r,(n+1)r;\R^d_+)$;

\item the solution is unique in the space $\cc(-r,(n+1)r;\R^d_+)$;

\item there exists a solution in the space $\cc(-r,(n+1)r;\R^d_+)$.
\end{enumerate}

\textit{Step} 1: \textit{If $x^{n+1}$ is a
solution of $\mathrm{(H_{n+1})}$ in the space $\cc(-r,(n+1)r;\R^d_+),$
then $x^{n+1}\in W_0^{\al,\infty}(-r,(n+1)r;\R^d_+)$.}  We can
write
\begin{eqnarray}\label{cotax}
 &&\|x^{n+1}\|_{\al,\infty(-r,(n+1)r)}\nonumber\\
 &&\quad=\sup
_{t\in[-r,(n+1)r]} \biggl(|x^{n+1}(t)|
 +\int_{-r}^t\frac{
|x^{n+1}(t)-x^{n+1}(s)|}{(t-s)^{\al+1}}\,\mathrm{d}s \biggr)\nonumber\\
 &&\quad\leq\sup_{t\in[-r,nr]}\biggl(
|x^{n}(t)|+\int_{-r}^t\frac{|x^{n}(t)-x^{n}(s)
|}{(t-s)^{\al+1}}\,\mathrm{d}s\biggr)\nonumber\\
 &&\qquad{}+\sup_{t\in[nr,(n+1)r]}
\biggl(|x^{n+1}(t)|+\int_{-r}^{nr}\frac{
|x^{n+1}(t)-x^{n}(s)|}{(t-s)^{\al+1}}\,\mathrm{d}s\nonumber\\
&&\phantom{\qquad{}+\sup_{t\in[nr,(n+1)r]}
\biggl(}{}+\int_{nr}^t\frac{
|x^{n+1}(t)-x^{n+1}(s)|}{(t-s)^{\al+1}}\,\mathrm{d}s\biggr)\\
&&\quad\leq\|x^{n}\|_{\al,\infty(-r,nr)}+\|x^{n+1}
\|_{\al,\infty(nr,(n+1)r)}\nonumber\\
&&\qquad{}+\sup_{t\in[nr,(n+1)r]}\int
_{-r}^{nr}\frac{|x^{n+1}(t)-x^n(nr)+x^n(nr)-x^{n}(s)
|}{(t-s)^{\al+1}}\,\mathrm{d}s\nonumber\\
&&\quad\leq2\|x^{n}\|_{\al,\infty(-r,nr)}+\|
x^{n+1}\|_{\al,\infty(nr,(n+1)r)}
 +\sup_{t\in[nr,(n+1)r]}\frac
{|x^{n+1}(t)-x^{n}(nr)|}{\al(t-nr)^{\al}}\nonumber\\
&&\quad =
2\|x^{n}\|_{\al,\infty(-r,nr)}+A_1+A_2,\nonumber
\end{eqnarray}
where
\begin{eqnarray*}
A_1&=& \|x^{n+1}\|_{\al,\infty(nr,(n+1)r)},\\
A_2&=&
\sup_{t\in[nr,(n+1)r]}\frac{|x^{n+1}(t)-x^{n}(nr)|}{\al
(t-nr)^{\al}}.
\end{eqnarray*}
%
From our hypothesis, it is clear that
$\|x^{n}\|_{\al,\infty(-r,nr)} < \infty$. So, to complete
this step of the proof, it suffices to check that $A_1 < \infty$ and
$A_2 <\infty$.

We begin with the study of $A_1$. Clearly,
\begin{eqnarray}\label{desa1}
A_1&\leq&
|\eta(0)|+\bigl\|F^{(b)}(x^{n+1})\bigr\|_{\al,\infty
(nr,(n+1)r)}\nonumber\\ [-6pt]\\ [-6pt]
&&{}+\bigl\|G^{(\si)}_r(x^n)\bigr\|_{\al
,\infty(nr,(n+1)r)}+\|y_{n+1}(\cdot)\|_{\al,\infty
(nr,(n+1)r)}.\nonumber
\end{eqnarray}
One of the keys to our proof is the study of the behavior of $y$. We
note that from its definition, it is clear that if $y^i$ is increasing
at $t$ (i.e., $y^i (t) > y^i (t - \varepsilon)$ for $\varepsilon\le
\varepsilon_0$ small enough), then $y^i(t)= - z^i (t)\geq0$. Moreover,
by construction, $y^i(s)\geq-z^i(s)$ for any $s$. So, if $y^i$ is
increasing at $t$, then for all $s<t,$
\[
|y^i(t)-y^i(s)|=y^i(t)-y^i(s)\leq-z^i(t)+ z^i(s)\leq|z^i(t)- z^i(s)|.
\]
For $t\in(nr,(n+1)r)$ and $i\in\{1,\ldots,d\}$, set
\[
t_0^i=\inf\{u; y^i(u)=y^i(t)\}\vee nr.
\]
Since $y^i$ is
increasing, we note that $y^i(s)=y^i(t_0^i)$ for all $s\in[t_0^i,t]$.
Then,
\begin{eqnarray}\label{cota1}
\int_{nr}^t
\frac{|y_{n+1}^i(t)-y_{n+1}^i(s)|}{(t-s)^{\al
+1}}\,\mathrm{d}s&=&\int_{nr}^{t_0^i}\frac{
|y_{n+1}^i(t_0^i)-y_{n+1}^i(s)|}{(t-s)^{\al+1}}\,\mathrm{d}s\nonumber\\
&\leq&\int_{nr}^{t_0^i}
\frac{|y_{n+1}^i(t_0^i)-y_{n+1}^i(s)|}{(t_0^i-s)^{\al
+1}}\,\mathrm{d}s\\
&\leq&\int_{nr}^{t_0^i}
\frac{|z^i_{n+1}(t_0^i)-z^i_{n+1}(s)|}{(t_0^i-s)^{\al+1}}\,\mathrm{d}s.\nonumber
\end{eqnarray}
On the other hand, we have
\begin{equation}\label{cota2}
|y_{n+1}^i(t)|= |y_{n+1}^i(t_0^i)|\leq\sup_{0
\le s \le t_0^i} |z^i_{n+1}(s)|.
\end{equation}
So, combining (\ref{cota1}) and (\ref{cota2}), we have that
\begin{equation}\label{cota3}
\|y_{n+1}\|_{\al,\infty(nr,(n+1)r)}\leq d
\bigl(\|z_{n+1}\|_{\al,\infty(nr,(n+1)r)}+\|
z_{n}\|_{\infty(0,nr)}\bigr),
\end{equation}
where we can use the bound
\begin{eqnarray}\label{cota4}
&&\|z_{n+1}\|_{\al,\infty(nr,(n+1)r)}\nonumber\\ [-6pt]\\ [-6pt]
&&\quad\leq
|\eta(0)|+\bigl\|F^{(b)}(x^{n+1})\bigr\|_{\al,\infty
(nr,(n+1)r)}+\bigl\|G_r^{(\si)}(x^n)\bigr\|_{\al
,\infty(nr,(n+1)r)}.\nonumber
\end{eqnarray}
Now combining (\ref{desa1}), (\ref{cota3}) and (\ref{cota4}), we
obtain that
\begin{eqnarray}\label{num2}
&&A_1 \le(d+1) \bigl( |\eta(0)|+\bigl\|F^{(b)}(x^{n+1})\bigr\|
_{\al,\infty(nr,(n+1)r)}\nonumber\\ [-6pt]\\ [-6pt]
&&\hspace*{-2.2pt}\hphantom{A_1 \le(d+1) \bigl(}{} + \bigl\|G_r^{(\si)}(x^n)\bigr\|_{\al,\infty(nr,(n+1)r)}\bigr) +
d\|z_{n}\|_{\infty(0,nr)}.\nonumber
\end{eqnarray}
From our hypothesis and Propositions \ref{Prop2.2} and \ref{Prop2.4},
it is easy to obtain that $ \|z_{n}\|_{\infty(0,nr)} <
\infty$. It therefore only remains to check the norms of the Lebesgue
and Riemann--Stieltjes integrals.

On the one hand,
%
%
\begin{eqnarray}\label{num3}
 &&\bigl\|F^{(b)}(x^{n+1})\bigr\|_{\al,\infty
(nr,(n+1)r)}\nonumber\\
&&\quad\leq\sup_{t\in[nr,(n+1)r]}\biggl(\int_0^t
|b(s,x^{n+1})|\,\mathrm{d}s
 +\int_{nr}^t\frac{\int_s^t
|b(u,x^{n+1})|\,\mathrm{d}u }{(t-s)^{\al+1}}\,\mathrm{d}s\biggr)\nonumber\\
&&\quad \leq\sup_{t\in[nr,(n+1)r]}\biggl(\int
_0^t\biggl(L_0\sup_{-r\leq u\leq s} |x^{n+1}(u)|+ b_0(s)
\biggr)\,\mathrm{d}s\\
&&\phantom{\quad \leq\sup_{t\in[nr,(n+1)r]}\biggl(}{} +\int_{nr}^t\frac{L_0\int
_s^t(\sup_{-r\leq v\leq u} |x^{n+1}(v)|+ b_0(u)
)\,\mathrm{d}u}{(t-s)^{\al+1}}\,\mathrm{d}s\biggr)\nonumber \\
 &&\quad\leq L_0\biggl(T+\frac{r^{1-\al}}{1-\al
}\biggr)\|x^{n+1}\|_{\infty(-r,(n+1)r)}
+\biggl(T^{1-\alpha}+\frac{r^{1-2\al}}{1-2\al}\biggr)\|
b_0\|_{L^{1/\al}}.\nonumber
\end{eqnarray}
 On the other hand, to study the Young integral, we will use
Proposition \ref{Prop2.4} and the fact that $x^n\in
W_0^{\al,\infty}(-r,nr;\R^d_+)$. For any $\la\geq1,$
\begin{eqnarray}\label{num4}
\bigl\|G_r^{(\si)}(x^n)\bigr\|_{\al,\infty(nr,(n+1)r)} & \leq&
\mathrm{e}^{\la(n+1)r}\bigl\|G_r^{(\si)}(x^n)\bigr\|_{\al,\la
(nr,(n+1)r)}\nonumber\\ [-8pt]\\ [-8pt]
&\leq&\frac{\Lambda_\al(g)d^{(4)}}{\la^{1-2\al}}\mathrm{e}^{\la(n+1)r}
\bigl(1+\|x^n\|_{\al,\la(-r,nr)}\bigr).\nonumber
\end{eqnarray}
%
So, combining (\ref{num2}), (\ref{num3}) and (\ref{num4}), and using
the facts that $\|x^{n}\|_{\alpha,\lambda(-r,nr)}<\infty$
and $\|x^{n+1}\|_{\infty(-r,(n+1)r)}<\infty$ for $\la
\geq
1$, we get that $A_1<\infty$.

We now deal with the
term $A_2$. We can write the decomposition
\begin{eqnarray}\label{novau}
\frac{|x^{n+1}(t)-x^n(nr)|}{(t-nr)^\al}&\leq&
\frac{|\int_{nr}^t
b(s,x^{n+1})\,\mathrm{d}s|}{(t-nr)^\al}+\frac{|\int_{nr}^t\si
(s,x^n(s-r))\,\mathrm{d}g_s|}{(t-nr)^\al}\nonumber\\ [-8pt]\\ [-8pt]
&&{}+\frac{
|y_{n+1}(t)-y_{n+1}(nr)|}{(t-nr)^\al}.\nonumber
\end{eqnarray}
 Using the same arguments as in (\ref{cota1}), we get
%
\begin{eqnarray}\label{cota7a}
 &&\sup_{t\in[nr,(n+1)r]} \frac{
|y_{n+1}(t)-y_{n+1}(nr)|}{(t-nr)^\al}\nonumber\\
&&\quad\leq d\sup_{t\in
[nr,(n+1)r]} \frac{|z_{n+1}(t)-z_{n+1}(nr)|}{(t-nr)^\al}\\
 &&\quad\leq d\biggl( \sup_{t\in[nr,(n+1)r]} \frac
{|\int_{nr}^t b(s,x^{n+1})\,\mathrm{d}s|}{(t-nr)^\al}
 +\sup_{t\in[nr,(n+1)r]}\frac
{|\int_{nr}^t
\si(s,x^{n}(s-r))\,\mathrm{d}g_s|}{(t-nr)^\al}\biggr).\nonumber
\end{eqnarray}
 Moreover, from Propositions \ref{Prop2.2} and \ref{Prop2.4},
we obtain the estimates
\begin{eqnarray}\label{cota7b}
 \sup_{t\in[nr,(n+1)r]} \frac{|\int_{nr}^t
b(s,x^{n+1})\,\mathrm{d}s|}{(t-nr)^\al}
&\leq&\bigl\|F^{(b)}(x^{n+1})\bigr\|_{1-\al(nr,(n+1)r)}r^{1-2\al}\nonumber\\ [-6pt]\\ [-6pt]
 &\leq&
\frac{d^{(1)}}{r^{2\al-1}}\bigl(1+\|x^{n+1}\|_{\infty
(-r,(n+1)r)}\bigr),\nonumber\\\label{cota7c}
 \sup_{t\in[nr,(n+1)r]}\frac{|\int_{nr}^t
\si(s,x^{n}(s-r))\,\mathrm{d}g_s|}{(t-nr)^\al}&\leq&
\bigl\|G^{(\si)}(x^n)\bigr\|_{1-\al(nr,(n+1)r)} r^{1-2\al}\nonumber\\[-6pt]\\ [-6pt]
 &\leq&
\frac{d^{(3)}\Lambda_\al(g)}{r^{2\al-1}}\bigl(1+\|x^n
\|_{\al,\infty(-r,nr)}\bigr).\nonumber
\end{eqnarray}
So, using the facts that $\|x^n\|_{\al,\infty(-r,nr)} <
\infty$ and $\|x^{n+1}\|_{\infty(-r,(n+1)r)}<\infty,$ and
combining~(\ref{novau}), (\ref{cota7a}), (\ref{cota7b}) and
(\ref{cota7c}), we obtain that $A_2 < \infty$.

 The
proof of the first step is now complete.\vspace*{6pt}

\textit{Step} 2: \textit{Uniqueness of the solution in the space
$\cc(-r,(n+1)r;\R^d_+)$.}

Let $x$ and $x'$ be two solutions of (\ref{equationn}) in the space
$\cc(-r,(n+1)r;\R^d_+)$ and choose $N$ large enough so that
$\|x\|_{\infty(-r,(n+1)r)}\leq N$ and
$\|x'\|_{\infty(-r,(n+1)r)}\leq N$.

For
any $t\in[0,(n+1)r]$,
\[
\sup_{s\in[0,t]} |x(s)-x'(s)| \le
\sup_{s\in[0,t]}|z(s)-z'(s)|+ \sup_{s\in[0,t]}
|y(s)-y'(s)|.
\]
Moreover, using Lemma \ref{le2}, we have
\[
\sup_{s\in[0,t]} |y(s)-y'(s)|\leq K_l \sup_{s\in[0,t]}
|z(t)-z'(t)|.
\]
So, combining the last two inequalities, we
get that
\begin{eqnarray*}
\sup_{s\in[0,t]} |x(s)-x'(s)|&\le& (1+K_l)
\sup_{s\in[0,t]} |z(s)-z'(s)| \\
&\leq& (1+K_l)
\sup_{s\in[0,t]} \biggl|\int_0^s
\bigl(b(u,x)-b(u,x')\bigr)\,\mathrm{d}u\biggr|\\
&\leq& (1+K_l)
L_N\sup_{s\in[0,t]}\biggl|\int_0^s\sup_{0\leq v\leq
u}|x(v)-x'(v)|\,\mathrm{d}u\biggr|\\
&\leq& L_N(1+K_l)\int_0^t\sup_{v\in[0,
u]}|x(v)-x'(v)|\,\mathrm{d}u.
\end{eqnarray*}
Now applying Gronwall's inequality, we have that for all
$t\in[0,(n+1)r],$
\[
\sup_{s\in[0,t]}|x(s)-x'(s)|= 0.
\]
So,
\[
\|x-x'\|_{\infty(-r,(n+1)r)}= 0
\]
and the uniqueness has been proven.\vspace*{6pt}

\textit{Step} 3: \textit{Existence of a solution in
$\cc(-r,(n+1)r;\R^d_+)$.}

In the space $\cc(-r,(n+1)r;\R^d_+),$ we can deal with the reflection
term using the Skorokhod mapping. Nevertheless, since the coefficient
$b$ is only locally Lipschitz, we will need to use a fixed point
argument in $\cc(-r,(n+1)r;\R^d_+)$ based on Lemma \ref{puntfix}.

Let us consider the operator
\[
\cl\dvtx \cc\bigl(-r,(n+1)r;\R^d_+\bigr)\rightarrow\cc\bigl(-r,(n+1)r;\R^d_+\bigr)
\]
which is such that
\begin{eqnarray*}
\cl(u)(t)&=&\eta(0)+\int_0^t
b(s,u)\,\mathrm{d}s+\int_0^t\si\bigl(s,x^n(s-r)\bigr)\,\mathrm{d}g_s+y_{n+1,u}(t),\qquad
t\in[0,(n+1)r],\\
\cl(u)(t)&=&\eta(t),\qquad t\in[-r,0],
\end{eqnarray*}
where $x^n$ is the solution obtained in $\mathrm{(H_{n})}$ and if
\[
z_{n+1,u}(t)=\eta(0)+\int_0^t b(s,u)\,\mathrm{d}s+ \int_0^t \si\bigl(s,x^n(s-r)\bigr)\,\mathrm{d}g_s,
\]
then ${
y_{n+1,u}^i(t)=\max_{s\in[0,t]}(z_{n+1,u}^i(s))^-}$ for all
$i=1,\ldots, d$.

Note that
$\cl$ is well defined. Moreover, if $u=\cl(u)$, then $u(t)=x^n(t)$
for any $t \in(-r,nr)$. We will use the notation $u^*=\cl(u)$.

We next introduce a new norm in the space $\cc(-r,(n+1)r;\R^d_+)$: for
any $\la\geq1,$
\[
\|f\|_{\infty,\la(-r,(n+1)r)}:= \sup_{t\in[-r,(n+1)r]}
\mathrm{e}^{-\la t}|f(t)|.
\]
These norms are equivalent to $\|f\|_{\infty
(-r,(n+1)r)}$.

We now check that we can apply Lemma \ref{puntfix}. Note first that
\begin{eqnarray}\label{cotaE1}
&&\hspace*{-30pt}\|u^*\|_{\infty,\la(-r,(n+1)r)}\nonumber\\
&&\hspace*{-30pt}\quad\leq \|\eta\|_{\infty,\la(-r,0)}+|\eta(0)|+\sup
_{t\in[0,(n+1)r]}\mathrm{e}^{-\la t}\biggl|\int_0^t b(s,u)\,\mathrm{d}s\biggr|\nonumber\\
&&\hspace*{-30pt}\qquad{}+\sup_{t\in[0,(n+1)r]}\mathrm{e}^{-\la t}\biggl|\int_0^t \si
\bigl(s,x^n(s-r)\bigr)\,\mathrm{d}g_s\biggr|+d\sup_{t\in[0,(n+1)r]}\mathrm{e}^{-\la t}|z_{n+1,u}(t)|\\
&&\hspace*{-30pt}\quad\leq
\|\eta\|_{\infty,\la(-r,0)}+(d+1)|\eta(0)|+
(d+1) \sup_{t\in[0,(n+1)r]}\mathrm{e}^{-\la t}\biggl|\int_0^t b(s,u)\,\mathrm{d}s
\biggr|\nonumber\\
&&\hspace*{-30pt}\qquad{}+(d+1) \sup_{t\in[0,(n+1)r]}\mathrm{e}^{-\la t}\biggl|\int_0^t
\si\bigl(s,x^n(s-r)\bigr)\,\mathrm{d}g_s\biggr|,\nonumber
\end{eqnarray}
where we have used computations similar to those in (\ref{cota2}).
Indeed, for fixed $t,$ let $t_1:=\inf\{u; y^i(u)=y^i(t)\}$.
Then,
\[
\mathrm{e}^{-\lambda t } \vert y_{n+1,u}^i (t) \vert\le \mathrm{e}^{-\lambda t_1 } \vert
y_{n+1,u}^i (t_1) \vert \le \mathrm{e}^{-\lambda t_1 } \vert z_{n+1,u}^i (t_1)
\vert,
\]
and taking suprema, we have
\begin{equation}\label{afesup}
\sup_{t \in[0,(n+1)r]} \mathrm{e}^{-\lambda t } \vert y_{n+1,u} (t) \vert \le
d \sup_{t \in[0,(n+1)r]} \mathrm{e}^{-\lambda t } \vert z_{n+1,u} (t)
\vert.
\end{equation}
Moreover, we have
\begin{eqnarray}\label{cotaE2}
&&\sup_{t\in[0,(n+1)r]} \mathrm{e}^{-\la t}\biggl|\int_0^t b(s,u)\,\mathrm{d}s\biggr|\nonumber\\
&&\quad\leq L_0
\sup_{t\in[0,(n+1)r]} \mathrm{e}^{-\la t}\int_0^t\sup_{-r\leq
v\leq s} |u(v)|\,\mathrm{d}s+ \sup_{t\in[0,(n+1)r]} \mathrm{e}^{-\la t}\biggl|\int_0^t
b_0(s)\,\mathrm{d}s\biggr|\\
&&\quad\leq\frac{L_0}{\la}
\|u\|_{\infty,\la(-r,(n+1)r)}+\frac{C_\al}{\la^{1-\al
}}\|b_0\|_{L^{1/\al}}\nonumber
\end{eqnarray}
and, from Proposition \ref{Prop2.4},
\begin{equation}\label{cotaE3}
\sup_{t\in[0,(n+1)r]} \mathrm{e}^{-\la t}\biggl|\int_0^t
\si\bigl(s,x^n(s-r)\bigr)\,\mathrm{d}g_s\biggr| \leq\frac{\Lambda_\al(g)
d^{(4)}}{\la^{1-2\al}}\bigl(1+\|x^n\|_{\al,\la
(-r,nr)}\bigr).
\end{equation}
So, combining (\ref{cotaE1}), (\ref{cotaE2}) and (\ref{cotaE3}), we
have
\[
\|u^*\|_{\infty,\la(-r,(n+1)r)}\leq
M_1(\la)+M_2(\la)\|u\|_{\infty,\la(-r,(n+1)r)},
\]
where
\begin{eqnarray*}
M_1(\la)&=&\|\eta\|_{\infty,\la(-r,0))}+(d+1)|\eta
(0)|+\frac{(d+1)C_\al}{\la^{1-\al}}\|b_0\|_{L^{1/\al}}\\
&&{}+(d+1)\frac{\Lambda_\al(g)d^{(4)}}{\la^{1-2\al}}
\bigl(1+\|x^n\|_{\al,\la(-r,nr)}\bigr),\\
M_2(\la)&=& (d+1)L_0\frac{1}{\la}.
\end{eqnarray*}
Choose $\la=\la_0$ large enough so that
$M_2(\la_0) \leq\frac{1}{2}$ . Then, if
${ \|u\|_{\infty,\la_0(-r,(n+1)r)}\leq
2M_1(\la_0)}$, we have
\[
{ \|u^*\|_{\infty,\la_0(-r,(n+1)r)}\leq
2M_1(\la_0)}
\]
and so ${ \cl(B_0)\subseteq B_0}$, where
\[
B_0=\bigl\{u\in\cc\bigl(-r,(n+1)r;\R^d_+\bigr); \|u\|_{\infty
,\la_0(-r,(n+1)r)}\leq2M_1(\la_0)\bigr\}.
\]
The first hypothesis in Lemma \ref{puntfix} is thus satisfied with the
metric $\rho_0$ associated with the norm
$\|\cdot\|_{\infty,\la_0(-r,(n+1)r)}$.\vadjust{\goodbreak}

To complete the proof, we only need to find a metric $\rho_1$
satisfying the second hypothesis of Lemma \ref{puntfix}.

 Note first that if $u\in B_0$, then
$\|u\|_{\infty(-r,(n+1)r)}\leq
2\mathrm{e}^{\la_0(n+1)r}M_1(\la_0):=N_0$. Consider $u,u'\in B_0$ and $\la
\geq
1$. We then have\vspace*{-1pt}
\begin{eqnarray*}
 \|\cl(u)-\cl(u')\|_{\infty,\la
(-r,(n+1)r)}&\leq&\sup_{t\in[0,(n+1)r]}\mathrm{e}^{-\la t}
|z_{n+1,u}(t)-z'_{n+1,u}(t)|\\[-2pt]
 &&{}+\sup_{t\in
[0,(n+1)r]}\mathrm{e}^{-\la
t}|y_{n+1,u}(t)-y'_{n+1,u}(t)|.\vspace*{-1pt}
\end{eqnarray*}
From Lemma \ref{le2}, note that given $t \in[0,(n+1)r],$ there exists
some $t_2 \le t$ such that\vspace*{-1pt}
\[
|y_{n+1,u}(t)-y'_{n+1,u}(t)| \le K_l
|z_{n+1,u}(t_2)-z'_{n+1,u}(t_2)| .\vspace*{-1pt}
\]
So,\vspace*{-1pt}
\[
\mathrm{e}^{-\lambda t} |y_{n+1,u}(t)-y'_{n+1,u}(t)| \le K_l
\mathrm{e}^{-\lambda t_2} |z_{n+1,u}(t_2)-z'_{n+1,u}(t_2)|\vspace*{-1pt}
\]
and it
easily follows that\vspace*{-1pt}
\[
\sup_{t\in[0,(n+1)r]}\mathrm{e}^{-\la t}
|y_{n+1,u}(t)-y'_{n+1,u}(t)|
\le K_l \sup_{t\in[0,(n+1)r]}\mathrm{e}^{-\la
t}|z_{n+1,u}(t)-z'_{n+1,u}(t)|.\vspace*{-1pt}
\]
Then,\vspace*{-1pt}
\begin{eqnarray*}
&& \|\cl(u)-\cl(u')\|_{\infty,\la
(-r,(n+1)r)}\\[-2pt]
&&\quad\leq(1+K_l) \sup_{t\in[0,(n+1)r]}\mathrm{e}^{-\la t}\int
_0^t|b(s,u)-b(s,u')|\,\mathrm{d}s\\[-2pt]
 &&\quad\leq L_{N_0} (1+K_l) \sup_{t\in
[0,(n+1)r]}\mathrm{e}^{-\la t}\int_0^t\sup_{0\leq v\leq s}
|u(v)-u'(v)|\,\mathrm{d}s\\[-2pt]
 &&\quad\leq L_{N_0}(1+K_l)\sup_{t\in
[0,(n+1)r]}\int_0^t \mathrm{e}^{-\la( t-s)} \mathrm{e}^{-\la s} \sup_{-r\leq v\leq
s}|u(v)-u'(v)|\,\mathrm{d}s\\[-2pt]
 &&\quad\leq L_{N_0}(1+K_l
)\frac{1}{\la}\|u-u'\|_{\infty,\la(-r,(n+1)r)} .\vspace*{-1pt}
\end{eqnarray*}
So, if we choose $\la=\la_1$ such that
$ \frac{L_{N_0}(1+K_l)}{\la}\leq\frac{1}{2}$, then the
second hypothesis is satisfied for the metric $\rho_1$ associated with
the norm $\|\cdot\|_{\infty,\la_1(-r,(n+1)r)}$ and
$a=\frac{L_{N_0}(1+K_l)}{\la_1}$.\vspace*{-2pt}
\end{pf}

We now check that the solution is $(1-\alpha)$-H\"{o}lder
continuous.\vspace*{-2pt}

\begin{proposition}\label{teohol}
Assume that b and $\sigma$ satisfy hypothesis \textup{(H1)}
and \textup{(H2)}, respectively, with $\rho=1/\alpha$ and $0
< \alpha<\min\{\frac12, \beta\}$. The solution $x$ of equation
(\ref{det}) then belongs to $C^{1-\al}(0,T;\R^d)$ with\vspace*{-1pt}
\[
\Vert x \Vert_{1-\alpha(0,T)} \le d^{(6)}\bigl(1+ \Delta_\alpha(g)\bigr) \bigl( 1
+ \Vert x \Vert_{\alpha,\infty(-r,T)}\bigr),\vspace*{-1pt}
\]
where $d^{(6)}$ is a positive constant independent of $f$ and
$g$.\vadjust{\goodbreak}
\end{proposition}
\begin{pf} Note that
\[
\Vert x \Vert_{1-\alpha(0,T)} \le\Vert z \Vert_{1-\alpha(0,T)} +
\Vert y \Vert_{1-\alpha(0,T)}.
\]
For fixed $t \in[0,T]$, set $t_*= \inf\{ u \le t; y^i (u) = y^i (t)
\}.$ Then, $y^i $ is increasing in $t_*$ and it is easy to check that
$\vert y^i (t_*) - y^i(s) \vert\le\vert z^i (t_*) - z^i(s) \vert$ for
all $s \in(0,t_*)$. For all $s \in(0,t_*),$ it thus holds that
\[
\frac{\vert y^i (t) - y^i(s) \vert}{(t-s)^{1-\alpha}} \le\frac
{\vert
z^i (t_*) - z^i(s) \vert}{(t_*-s)^{1-\alpha}}
\]
and it then follows
easily that $\Vert y \Vert_{1-\alpha(0,T)} \le d \Vert z
\Vert_{1-\alpha(0,T)}$. So,
\begin{eqnarray*}
\Vert x \Vert_{1-\alpha(0,T)} & \le& (d+1) \Vert z \Vert_{1-\alpha
(0,T)} \\[-2pt]
&\le& (d+1) \bigl( \vert\eta(0) \vert+ \bigl\Vert F^{(b)}(x)
\bigr\Vert_{1-\alpha(0,T)} + \bigl\Vert G_r^{(\sigma)}(x)
\bigr\Vert_{1-\alpha(0,T)} \bigr).
\end{eqnarray*}
Using Propositions \ref{Prop2.2} and \ref{Prop2.4}, we easily complete
the proof.
\end{pf}

We will now give an upper bound for the norm of
the solution. Recall the definition of~$\varphi(\gamma,\al)$:
\[
\varphi(\gamma,\al)=
\cases{
2\al,&\quad$\gamma=1$,\vspace*{1pt}\cr
\displaystyle>1+\frac{2\al-1}{\gamma},&\quad$\displaystyle\frac{1-2\al}{1-\al}\leq\gamma<
1$,\cr
\displaystyle\al,&\quad $\displaystyle0\leq\gamma< \frac{1-2\al}{1-\al}$.
}
\]

\begin{lemma}
Assume \textup{(H1)}, \textup{(H2)} and
\textup{(H3)}. The unique solution of equation (\ref{det})
then satisfies
\[
\|x\|_{\al,\infty(r)}\leq d_\al^{(3)}\bigl(\|
\eta\|_{\al,\infty(-r,0)}+\Lambda_\al(g)+1\bigr)\exp
\bigl(T\bigl(d^{(1)}_\al+d^{(2)}_\al\Lambda_\al(g)^{1/(1-\varphi
(\gamma,\al))}\bigr)\bigr).
\]
\end{lemma}
\begin{pf} First, we need to obtain an upper bound for
$\|x\|_{\al,\la(r)}$. We begin with the estimates
\begin{eqnarray}\label{dd1}
\|x\|_{\al,\la(r)}&\leq& \|\eta\|_{\al
,\la(-r,0)}+\sup_{t\in[0,T]} \mathrm{e}^{-\la t}\biggl(|x(t)|+\int_{-r}^t
\frac{|x(t)-x(s)|}{(t-s)^{\al+1}}\,\mathrm{d}s\biggr)\nonumber\\[-2pt]
&\leq& \|\eta\|_{\al,\la(-r,0)}+\sup_{t\in[0,T]}
\mathrm{e}^{-\la t}\biggl(|x(t)|+\int_{-r}^0 \frac{|x(t)-\eta(s)|}{(t-s)^{\al
+1}}\,\mathrm{d}s\nonumber\\ [-9pt]\\ [-9pt]
&&\hphantom{\|\eta\|_{\al,\la(-r,0)}+\sup_{t\in[0,T]}
\mathrm{e}^{-\la t}\biggl(}{}+\int_{0}^t \frac
{|x(t)-x(s)|}{(t-s)^{\al+1}}\,\mathrm{d}s\biggr)\nonumber\\[-2pt]
&\leq& \|\eta\|_{\al,\la(-r,0)} +
\|x\|_{\al,\la(0,T)}+\sup_{t\in[0,T]} \mathrm{e}^{-\la t}\int_{-r}^0
\frac{|x(t)-\eta(s)|}{(t-s)^{\al+1}}\,\mathrm{d}s.\nonumber
\end{eqnarray}
Moreover,
\begin{equation}\label{dd2}
\|x\|_{\al,\la(0,T)}\leq
\|z\|_{\al,\la(0,T)}+\|y\|_{\al,\la
(0,T)}\vadjust{\goodbreak}
\end{equation}
and using the same arguments as in (\ref{cota1}) and (\ref{cota2}), we
have,
\begin{equation}\label{dd3}
\|y\|_{\al,\la(0,T)}\leq d
\|z\|_{\al,\la(0,T)}.
\end{equation}
We also know that
\begin{equation}\label{dd4}
\|z\|_{\al,\la(0,T)}\leq|\eta(0)|+
\bigl\|F^{(b)}(x)\bigr\|_{\al,\la(0,T)}+\bigl\|G^{(\si
)}_r(x)\bigr\|_{\al,\la(0,T)}.
\end{equation}
So, combining (\ref{dd1})--(\ref{dd4}) and applying Propositions
\ref{Prop2.2} and \ref{Prop2.6}, we get that
\begin{eqnarray}\label{dd5}
\|x\|_{\al,\la(r)}&\leq&
\|\eta\|_{\al,\la(-r,0)}+(d+1) |\eta(0)|+
(d+1)d^{(2)}\biggl(
\frac{1}{\lambda^{1-2\alpha}} +
\frac{\|x\|_{\al,\la(r)}}{\lambda^{1-\alpha}}\biggr)
\nonumber\\ [-6pt]\\ [-6pt]
&&{}+\Lambda_\al(g)(d+1)d^{(5)}\biggl(1+\frac{\|x\|_{\al
,\la(r)}}{\la^{1-\varphi(\gamma,\al)}}\biggr)
+ B \nonumber
\end{eqnarray}
with
\[
B:=\sup_{t\in[0,T]} \mathrm{e}^{-\la t}\int_{-r}^0
\frac{|x(t)-\eta(s)|}{(t-s)^{\al+1}}\,\mathrm{d}s.
\]
 It remains to investigate $B$. We can decompose this term as
follows:
\begin{eqnarray}\label{dd6}
B&\leq& \sup_{t\in[0,T]} \mathrm{e}^{-\la t}\int_{-r}^0 \frac{|x(t)-\eta
(0)|}{(t-s)^{\al+1}}\,\mathrm{d}s+\sup_{t\in[0,T]} \mathrm{e}^{-\la t}\int_{-r}^0 \frac
{|\eta(0)-\eta(s)|}{(-s)^{\al+1}}\,\mathrm{d}s\nonumber\\[2pt]
&\leq& \frac{1}{\al}\sup_{t\in[0,T]} \frac{\mathrm{e}^{-\la t}}{t^\al
}|x(t)-\eta(0)|+\|\eta\|_{\al,\la(-r,0)}\\[2pt]
&\leq& \frac{1}{\al} (B_1+B_2+B_3)+
\|\eta\|_{\al,\la(-r,0)},\nonumber
\end{eqnarray}
where
\begin{eqnarray*}
B_1&=& \sup_{t\in[0,T]}\frac{\mathrm{e}^{-\la t}}{t^\al}\int_0^t
|b(s,x)|\,\mathrm{d}s,\\[2pt]
B_2&=& \sup_{t\in[0,T]}\biggl|\frac{\mathrm{e}^{-\la t}}{t^\al}\int_0^t \si
\bigl(s,x(s-r)\bigr)\,\mathrm{d}g_s\biggr|,\\[2pt]
B_3&=& \sup_{t\in[0,T]} \frac{\mathrm{e}^{-\la t}}{t^\al}|y(t)|.
\end{eqnarray*}
 Using the same arguments as in (\ref{afesup}), we get that
\[
B_3=\sup_{t\in[0,T]} \frac{\mathrm{e}^{-\la t}}{t^\al}|y(t)|\leq d \sup
_{t\in[0,T]} \frac{\mathrm{e}^{-\la t}}{t^\al}|z(t)|\leq d(B_1+B_2).
\]
 We now consider $B_1$ and $B_2$. For $B_1,$ we can write
\begin{eqnarray*}
B_1&\leq& \sup_{t\in[0,T]} \frac{\mathrm{e}^{-\la t}}{t^\al}\int_0^t
\Bigl(L_0\sup_{-r\leq u\leq s} |x(u)|+b_0(s)\Bigr)\,\mathrm{d}s\\
&\leq& L_0\Bigl(\sup_{s\in[-r,T]} \mathrm{e}^{-\la s}|x(s)|\Bigr)\sup
_{t\in[0,T]}\int_0^t\frac{\mathrm{e}^{-\la(t-s)}}{(t-s)^\al}\,\mathrm{d}s+\sup_{t\in[0,T]} \frac{\mathrm{e}^{-\la t}}{t^{2\al-1}}\|b_0
\|_{L^{1/\al}}\\
&\leq& L_0\la^{\al-1}\Gamma(1-\al)\|x\|_{\al,\la(r)}
+C_\al\la^{2\al-1}\|b_0\|_{L^{1/\al}}.
\end{eqnarray*}
Next, we obtain a bound for $B_2$. We will use the hypothesis
(H3).
\begin{eqnarray*}
B_2 &\leq& \sup_{t\in[0,T]} \frac{\mathrm{e}^{-\la t}}{t^\al}\Lambda_\al
(g)\biggl(\int_0^t\frac{|\si(s,x(s-r))|}{s^\al}\,\mathrm{d}s\\
&&\phantom{\sup_{t\in[0,T]} \frac{\mathrm{e}^{-\la t}}{t^\al}\Lambda_\al
(g)\biggl(}{}+\al\int_0^t\int_0^s\frac{|\si(s,x(s-r))-\si
(y,x(y-r))|}{(s-y)^{\al+1}}\,\mathrm{d}y\,\mathrm{d}s\biggr)\\
&\leq& \sup_{t\in[0,T]} \frac{\mathrm{e}^{-\la t}}{t^\al}\Lambda_\al
(g)\biggl(K_0\int_0^t\frac{1+|x(s-r)|^\gamma}{s^\al}\,\mathrm{d}s\\
&&\hphantom{\sup_{t\in[0,T]} \frac{\mathrm{e}^{-\la t}}{t^\al}\Lambda_\al
(g)\biggl(}{}+\al M_0\int_0^t\int_0^s\biggl(\frac
{|x(s-r)-x(y-r)|}{(s-y)^{\al+1}}+\frac{1}{(s-y)^{\al+1-\beta
}}\biggr)\,\mathrm{d}y\,\mathrm{d}s\biggr)\\
&\leq& \sup_{t\in[0,T]} \Lambda_\al(g)\biggl(K_0\frac{t^{1-2\al
}}{1-\al}\mathrm{e}^{-\la t}+K_0\frac{\mathrm{e}^{-\la t}}{t^\al}\int_{-r}^{t-r}\frac
{|x(s)|^\gamma}{(s+r)^\al}\,\mathrm{d}s\\
&&\hphantom{\sup_{t\in[0,T]} \Lambda_\al(g)\biggl(}{}+\al M_0\|x\|_{\al,\la(r)}\int_{-r}^{t-r}
\frac{\mathrm{e}^{-\la(t-s)}}{(t-s)^\al}\,\mathrm{d}s +\frac{\al M_0t^{\beta-2\al
+1}\mathrm{e}^{-\la
t}}{(\beta-\al)(\beta-\al+1)}\biggr).
\end{eqnarray*}
Now, using the inequalities
\begin{eqnarray*}
\sup_{t\in[0,T]} t^\mu \mathrm{e}^{-\la t}&\leq&\biggl(\frac{\mu}{\la}
\biggr)^\mu
\mathrm{e}^{-\mu},\\
\int_{-r}^{t-r} \frac{\mathrm{e}^{-\la(t-s)}}{(t-s)^\al}\,\mathrm{d}s&=& \mathrm{e}^{-\la r}\int
_{0}^{t} \frac{\mathrm{e}^{-\la(t-u)}}{(t-u+r)^\al}\,\mathrm{d}u\\
&\leq& \mathrm{e}^{-\la r}\int_{0}^{t} \frac{\mathrm{e}^{-\la(t-u)}}{(t-u)^\al
}\,\mathrm{d}u\leq
\mathrm{e}^{-\la r}\la^{\al-1}\Gamma(1-\al),
\end{eqnarray*}
the H\"{o}lder inequality and the fact that $|f(s)|^\gamma\leq
|f(s)|+1,$ we get
\begin{eqnarray*}
\frac{\mathrm{e}^{-\la t}}{t^\al}\int_{-r}^{t-r}\frac{|x(s)|^\gamma
}{(s+r)^\al}\,\mathrm{d}s&\leq& \mathrm{e}^{-\la t}t^{\varphi(\gamma,\al)\gamma-2\al
+1-\gamma}\biggl(\int_{-r}^{t-r}\frac{|x(s)|}{(s+r)^{\varphi(\gamma
,\al)}}\,\mathrm{d}s\biggr)^\gamma\\
&\leq& C_{\al,\gamma,T}\mathrm{e}^{-\la t}\biggl(1+\int_{-r}^{t-r}\frac
{|x(s)|}{(s+r)^{\varphi(\gamma,\al)}}\,\mathrm{d}s\biggr)\\
&\leq& C_{\al,\gamma,T}\biggl(1+\|x\|_{\al,\la
(r)}\mathrm{e}^{-\la r}\int_{0}^{t}\frac{\mathrm{e}^{-\la(t-u)}}{u^{\varphi(\gamma
,\al)}}\,\mathrm{d}u\biggr)\\
&\leq& C_{\al,\gamma,T}\bigl(1+\|x\|_{\al,\la
(r)}\mathrm{e}^{-\la
r}\la^{\varphi(\gamma,\al)-1}\bigr),
\end{eqnarray*}
where we have used the fact that
$\varphi(\gamma,\al)\gamma-2\al+1-\gamma\geq0$. We thus obtain
\begin{eqnarray*}
B_2&\leq& \Lambda_\al(g)\biggl(\frac{K_0}{1-\al}\biggl(\frac
{1-2\al}{\mathrm{e}}\biggr)^{1-2\al}\la^{2\al-1}+\frac{\al M_0((\beta-2\al+1)\mathrm{e})^{\beta-2\al
+1}}{(\beta-\al)(\beta-\al+1)}\la^{2\al-1-\beta}\\
&&\hphantom{\Lambda_\al(g)\biggl(}{}+C_{\al,\gamma,T}+\|x\|_{\al,\la(r)} \mathrm{e}^{-\la
r}\bigl(\al M_0\Gamma(1-\al)\la^{\al-1}+K_0C_{\al,\gamma,T}\la
^{\varphi(\gamma,\al)-1}\bigr)\biggr)\\
&\leq&\Lambda_\al(g)C_{\al,\beta,\gamma}\bigl(1+\la^{2\al-1}+\mathrm{e}^{-\la
r}\la^{\varphi(\gamma,\al)-1}\|x\|_{\al,\la(r)}\bigr).
\end{eqnarray*}
%
Finally, combining (\ref{dd5}), (\ref{dd6}) and the estimates for $B,
B_1$ and $B_2$, we have
\[
\|x\|_{\al,\la(r)}\leq
M_1(\la)+M_2(\la)\|x\|_{\al,\la(r)}
\]
with
\begin{eqnarray*}
M_1(\la)&=& 2\|\eta\|_{\al,\la(-r,0)}+(d+1)
\biggl(|\eta(0)|+\Lambda_\al(g)d^{(5)}+C_{\al,\beta,\gamma}\\
&&\hphantom{2\|\eta\|_{\al,\la(-r,0)}+(d+1)
\bigl(}{}+\frac{C_{\al,\beta}}{\la^{1-2\al}}\bigl(d^{(2)}+\|
b_0\|_{L^{1/\al}}+\Lambda_\al(g)\bigr)\biggr),\\
M_2(\la)&=&
\frac{(d+1)C_\al}{\la^{1-\varphi(\gamma,\al)}}\bigl(\Lambda_\al
(g)\bigl(d^{(5)}+C_{\al,\beta,\gamma}\bigr)+L_0\Gamma(1-\al
)+d^{(2)}\bigr).
\end{eqnarray*}
Choosing $\la=\la_0$ large enough so that $M_2(\la_0)= \frac
{1}{2}$, we
then have
\[
\|x\|_{\al,\la_0(r)}\leq2M_1(\la_0).
\]
Set
\begin{eqnarray*}
\la_0&=& \bigl[2C_{\al,d}\bigl(d^{(2)}+\Lambda_\al
(g)\bigl(d^{(5)}+1\bigr)+L_0\bigr)\bigr]^{1/(1-\varphi(\gamma,\al
))}\\
&\leq& d_\al\bigl(2C_{\al,d}\bigl(d^{(2)}+L_0\bigr)
\bigr)^{1/(1-\varphi(\gamma,\al))}+\Lambda_\al(g)^{1/(1-\varphi(\gamma,\al))}d_\al
\bigl(2C_{\al,d}\bigl(1+d^{(5)}\bigr)\bigr)^{1/(1-\varphi
(\gamma,\al))}\\
&\leq&
d^{(1)}_\al+d^{(2)}_\al\Lambda_\al(g)^{1/(1-\varphi(\gamma
,\al))}
\end{eqnarray*}
with
\begin{eqnarray*}
d^{(1)}_\al&=& d_\al\bigl(2C_{\al,d}\bigl(d^{(2)}+L_0
\bigr)\bigr)^{1/(1-\varphi(\gamma,\al))},\\
d^{(2)}_\al&=&
d_\al\bigl(2C_{\al,d}\bigl(1+d^{(5)}\bigr)\bigr)^{
1/(1-\varphi(\gamma,\al))}.
\end{eqnarray*}
This implies that
\[
\|x\|_{\al,\infty(r)}\leq
\exp\bigl(T\bigl(d^{(1)}_\al+d^{(2)}_\al\Lambda_\al(g)^{1/(1-\varphi
(\gamma,\al))}\bigr)\bigr)2
M_1(\la_0)
\]
and the proof is then easily completed. Note that we can choose
$d^{(1)}_\al, d^{(2)}_\al$ which do not depend on $\beta$ or $\gamma$.
\end{pf}

\section{Stochastic integral equations}\label{sec5}

In this section, we
apply the deterministic results
in order to prove the main theorem of this paper.

 The stochastic integral appearing throughout this paper, $
\int_0^T u(s)\,\mathrm{d}W_s $, is a pathwise Riemann--Stieltjes integral and it
is well know that this integral exists if the process $u(s)$ has
H\"{o}lder continuous trajectories of order larger than $1-H$.

Set $\al\in(1-H, \frac12)$. For any $\delta\in(0,2)$, by Fernique's
theorem (see \cite{F}, Theorem 1.3.2), we have
\[
E( \exp(\Lambda_\al
(W)^\delta)) < \infty.
\]
Then, if $u=\{u_t, t \in[0,T]\}$ is a
stochastic process whose trajectories belong to the space $W_T^{\al,1}
(0,T)$, it follows that the Riemann--Stieltjes integral $\int_0^T
u(s)\,\mathrm{d}W_s$ exists and we have that
\[
\biggl| \int_0^T u(s)\,\mathrm{d}W_s \biggr| \le G \Vert u \Vert_{\al,1},
\]
where $G$ is a random variable with moments of all orders (see
\cite{NR}, Lemma 7.5). Moreover, if the trajectories of $u$ belong to
$W_0^{\al,\infty} (0,T)$, then the indefinite integral $\int_0^T u(s)\,\mathrm{d}W_s $ is H\"{o}lder continuous of order $1-\al$ and with trajectories
in $W_0^{\al,\infty} (0,T)$.

The existence and
uniqueness of a solution then follows from Theorem \ref{teoexi}. In
order to get the existence of a moment of any order, we need only note
that if $\al< (2-\gamma)/4,$ then $1/(1-\varphi(\gamma,\al)) < 2$ and
$ E( \exp(C \Lambda_\al(W)^{1/(1-\varphi(\gamma,\al))} )) <
\infty.$

\begin{appendix}\label{appendix}

\section*{Appendix}

In this appendix, we just give a fixed point theorem, well posed to our
problem, and recall a result with some properties of the solution of
the Skorokhod problem.

\begin{lemma} \label{puntfix}
Let $(X,\rho)$ be a complete metric space, and $\rho_0$ and $\rho_1$
two metrics on X which are equivalent to $\rho$. If $\cl\dvtx
X\rightarrow
X$ satisfies:
\begin{enumerate}
\item there exists some $r_0>0$, $x_0\in X$ such that if
$B_0=\{x\in X; \rho_0(x_0,x)\leq r_0\}$ then
$\cl(B_0)\subseteq B_0$;
\item there exists some $a\in(0,1)$ such that
$\rho_1(\cl(x),\cl(y))\leq a\rho_1(x,y)$ for all
$x,y\in
B_0$,
\end{enumerate}
then there exists some $x^*\in\cl(B_0)\subseteq X$ such that
$x^*=\cl(x^*)$.
\end{lemma}
\begin{pf} For all $n\in\N$, set
\[
x_{n+1}=\cl(x_n).
\]
Clearly, $x_n\in\cl(B_0)$ for all $n\in \N$. Moreover,
\[
\rho_1(x_{n+1}, x_n)=\rho_1(\cl(x_n),\cl(x_{n-1}))\leq
a\rho_1(x_n,x_{n-1})\leq\cdots\leq a^n\rho_1(x_1,x_0)
\]
and
\begin{eqnarray*}
\rho_1(x_{n+p},x_n)&\leq& \rho_1(x_{n+p},
x_{n+p-1})+\cdots+\rho_1(x_{n+1}, x_n)\\
&\leq&
a^n(a^{p-1}+\cdots+a+1)\rho_1(x_1, x_0) \leq
\frac{a^{n}}{1-a}\rho_1(x_1,x_0)\rightarrow0
\end{eqnarray*}
as $n\rightarrow\infty$.

Since $(X,\rho)$ is a complete metric space and $B_0$ is closed in X,
there exists some $x^*\in B_0$ such that $x_n\rightarrow x^*$.
Furthermore, from the second hypothesis of the lemma, we get that
\[
\rho_1(\cl(x_n),\cl(x^*))\leq a\rho_1(x_n,x^*).
\]
Since $\rho_1(x_n,x^*)\rightarrow0$, $\cl(x_n)\rightarrow\cl(x^*)$
and it follows that $x^*=\cl(x^*)$.
\end{pf}

\begin{lemma} \label{le2}
For each path $z\in\cc(\R_+,\R^d)$, there exists a unique
solution
$(x,y)$ to the Skorokhod problem for z. Thus, there exists a pair of
functions $ (\phi,\varphi)\dvtx
\cc_+(\R_+,\R^d)\rightarrow\break \cc_+(\R_+,\R^{2d})$ defined by
$(\phi(z),\varphi(z))=(x,y)$ such that the following holds:
 There exists a~constant $K_l>0$ such that for any
$z_1,z_2\in\cc_+(\R_+,\R^d),$ we have, for each $t\geq0$,
\begin{eqnarray*}
\|\phi(z_1)-\phi(z_2)\|_{\infty,[0,t]}&\leq& K_l
\|z_1-z_2\|_{\infty,[0,t]},\\
\|\varphi(z_1)-\varphi(z_2)\|_{\infty,[0,t]}&\leq&
K_l\|z_1-z_2\|_{\infty,[0,t]}.
\end{eqnarray*}
\end{lemma}

\begin{pf}
See \cite{KW}, Proposition A.0.1.
\end{pf}
\end{appendix}

\section*{Acknowledgement}
This work was partially supported by DGES Grant MTM09-07203 (both
authors).

\printhistory

\end{document}